%% file: main.tex
\documentclass[reqno]{siamart220329}
\usepackage[english]{babel}
\usepackage{lipsum}
\makeatletter

\newcommand{\authorfootnotes}{\renewcommand\thefootnote{\@fnsymbol\c@footnote}}%
\makeatother
\usepackage{etoolbox} %

\newcommand*\linenomathpatch[1]{%
  \cspreto{#1}{\linenomath}%
  \cspreto{#1*}{\linenomath}%
  \csappto{end#1}{\endlinenomath}%
  \csappto{end#1*}{\endlinenomath}%
}

\linenomathpatch{equation}
\linenomathpatch{gather}
\linenomathpatch{multline}
\linenomathpatch{align}
\linenomathpatch{alignat}
\linenomathpatch{flalign}

\usepackage{amsfonts}
\usepackage{amssymb}
\usepackage{centernot}
\usepackage{graphicx}
\usepackage{algorithm}
\usepackage{algorithmic}
\usepackage[algo2e]{algorithm2e}

\usepackage{fullpage}
\usepackage{lineno}

\usepackage{url}

\usepackage{color}

\DeclareMathOperator{\arginf}{arginf\;}
\DeclareMathOperator{\subto}{subject\,to\;}
\title{Learning solutions to some toy constrained optimization problems in infinite dimensional Hilbert spaces}

\begin{document}


\maketitle

 {\normalsize
 \centering
  \authorfootnotes
  Pinak Mandal\footnote{\thanks{Corresponding author: \texttt{pinak.mandal@icts.res.in}}}\textsuperscript{1,2}

  \textsuperscript{1} International Centre for Theoretical Sciences - TIFR, Bangalore 560089 India \par
   \textsuperscript{2}The University of Sydney, NSW 2050  Australia \par
 }

\begin{abstract}
In this work we present deep learning implementations of two popular theoretical constrained optimization algorithms in infinite dimensional Hilbert spaces, namely, the penalty and the augmented Lagrangian methods. We test these algorithms on some toy problems originating in either calculus of variations or physics. We demonstrate that both methods are able to produce decent approximations for the test problems and are comparable in terms of different errors produced. Leveraging the common occurrence of the Lagrange multiplier update rule being computationally less expensive than solving subproblems in the penalty method, we achieve significant speedups in cases when the output of the constraint function is itself a function.
\end{abstract}

\section{Introduction}
\label{sec-intro}
\input{var-tex/intro}

\section{Problem Statement and examples}
\label{sec-prob}
\input{var-tex/problem}

\section{Methodology}
\label{sec-method}
\input{var-tex/method}

\section{Results}
\label{sec-results}
\input{var-tex/results}

\section{Summary and future work}
\label{sec-conclusions}
\input{var-tex/conclusion}

\section*{Acknowledgements}
\input{var-tex/acknowledge}

\section{Appendix}
\label{sec-appendix}
\input{var-tex/appendix}

\bibliographystyle{siamplain}
\bibliography{ref}

\end{document}

%% file: var-tex/intro.tex
Pierre de Fermat authored many important works on his method of maxima and minima, out of which the last two were titled \textit{The analysis of refractions} and \textit{The synthesis of refractions}. These contained derivations of the law of refraction now commonly known as \textit{Snell's law}. In these papers Fermat states his intuition about the nature of physical laws as, \textit{"nature operates by means and ways that are easiest and
fastest"} \cite{goldstine2012history}. Even though the ancient Greeks had considered some classic problems in calculus of variations such as  isoperimetric problems \cite{pappus1933pappus}, Fermat's proclamations are one of the first instances where we encounter the notion that the laws of physics can often be stated in terms of optimization problems. This notion takes its final form as the principle of stationary action in modern physics appearing in nearly every subfield from classical mechanics, thermodynamics, relativity, quantum mechanics, string theory and everything in between \cite{coopersmith2017lazy}, \cite{rojo2018principle}. Since Fermat's time calculus of variations in general has found applications in most fields dealing with mathematical models, be it chemistry \cite{quapp2008chemical} or economics \cite{guzowska2015calculus} and its stochastic counterpart is useful in economics \cite{oksendal1997introduction} and mathematical finance \cite{malliavin2006stochastic}. 

Calculus of variations deals with finding functions as optimizer of functionals under constraints. Due to recent technical advancements in automatic differentiation and machine learning, it has become a popular paradigm to cast many engineering or basic science problems such as finding language models for Shakespearean text \cite{jhamtani2017shakespearizing}, learning generative models for natural images \cite{huang2018introduction}, solving partial differential equations \cite{blechschmidt2021three} etc as optimization problems and then solve them using well-established optimization algorithms like stochastic gradient descent \cite{ruder2016overview}, \cite{bottou2007tradeoffs}. This pattern very naturally yields itself to the \textit{function-finding} problems of calculus of variations. In finite dimensions, constrained optimization problems are routinely handled with \textit{penalty method}, \textit{augmented Lagrangian method} and their many variants \cite{jorge2006numerical}, \cite{bonnans2006numerical}, \cite{bertsekas1995athena}. In infinite dimensions or for function finding problems analogues of these algorithms have been discussed extensively in terms of theory \cite{ito2008lagrange}, \cite{kanzow2018augmented}, \cite{dussault2007penalty}, \cite{fiacco1969generalized}. But numerical implementation of these algorithms remain few and far between. This work aims to bridge the gap between the theory of infinite dimensional constrained optimization algorithms and their practical implementations using deep learning. Recently variational problems with essential boundary conditions have been explored by Huang et al \cite{huang2021augmented}. In this work we explore more general problems. Our goal is to evaluate our algorithms, rather than solving the specific problems we list here. We, therefore, apply them on some simple toy problems with known solutions. Our problems are either taken from the classics in calculus of variations or inspired by physics.

%% file: var-tex/problem.tex
In this work, we are interested in problems of the following form.
\begin{equation}
\begin{aligned}
    &\underset{u\in X}\arginf f(u)\\
    &\subto g(u)=0 
\end{aligned}   \label{eq:con-opt}
\end{equation}
where $f :X\to\mathbb R$ and $g: X\to W$ and $X, W$ are real Hilbert spaces. $X$, in particular, is an infinite dimensional Hilbert space whereas $W$ can be either finite or infinite dimensional. This ensures that problem~\eqref{eq:con-opt} is indeed an infinite-dimensional optimization problem. This setup allows us to encompass a fairly large class of problems with one or multiple constraints or even unconstrained problems if we set $g$ to be the zero function. 
To better familiarize ourselves with this setup let us first look at a few examples.
\subsection{The minimal surface problem} During the later half of the eighteenth century Lagrange in his correspondence with Euler delineated the foundations of calculus of variations and derived the famous Euler-Lagrange formula \cite{goldstine2012history}. One of the problems he considered during this time asks to find the surface of least area stretched across a given contour. Although Lagrange did not find any solutions other than the plane, Euler and Jean Baptiste Meusnier later showed that helicoid and catenoid are also valid solutions to the minimal surface problem \cite{meusnier1785memoire}. Since then the theory of minimal surfaces has seen multiple revivals with Schwarz's solution to the Björling problem \cite{darboux1896leccons}, the discovery of Costa's surface \cite{costa1984example} and has even found its way into mathematical physics through topics like positive energy theorem \cite{schoen1979proof}. The rich theory behind minimal surfaces allows them to be expressed in many different ways \cite{colding2011course}. Here we will work with a definition that closely resembles Lagrange's original formulation. Rather than describing the minimal surface problem in its full generality, we describe the specific problem we will solve below. We define $X$ to be an appropriate Sobolev space, $f$ to be an area functional and $g$ to be the boundary condition.
\begin{equation}
\begin{aligned}
    &\Omega=(0,1 )\times(-2\pi, 2\pi),\;X=W^{1, 2}(\Omega;\mathbb R),\; W=L^{2}((-2\pi,2\pi);\mathbb R)\\
    &f(u)=\int_{-2\pi}^{2\pi}\int_0^1\sqrt{\left[1+\left(\frac{\partial u}{\partial r}\right)^2\right]r^2+\left(\frac{\partial u}{\partial\theta}\right)^2}\;dr\,d\theta\\
&g(u): \theta \mapsto u(1, \theta) - \theta\label{eq:ms}
\end{aligned}    
\end{equation}
Here $W^{k, p}$ denotes the Sobolev space of function with $k$ $p$-integrable weak derivatives. Our question thus becomes, what is the surface of minimal area given it has a unit helix as its boundary? The solution $u^*$ gives us the minimal surface $(r, \theta, u^*(r,\theta))$. Note that even though we have used the standard area integral in polar coordinates, we are working beyond the standard domain of $\theta$ which is $[0, 2\pi)$. Therefore, when we visualize the solution to this problem using standard polar to Cartesian conversion we get a multivalued function or a helicoid with two full twists rather than just one, as seen in section~\ref{ssec:res-ms}.

\subsection{Geodesics on a surface}  Johann Bernoulli was interested in several problems in calculus of variations and investigated both curves of shortest length and time between two points \cite{struik1961lectures}, \cite{goldstine2012history}. The former type of curves are known as geodesics while the latter are known as brachistochrones. After having found the solution to the brachistochrone problem Bernoulli had challenged his contemporaries to come up with their own solutions (a practice that was not uncommon in the era) to which Newton (anonymously), Jacob Bernoulli, Leibniz and de L'Hôpital had responded with their own solutions. The aftermath of this challenge would eventually lead to the infamous calculus controversy between Leibniz and Newton \cite{palomo2021new}, \cite{goldstine2012history}. Even though the brachistochrone problem is one of the oldest problems to be posed in calculus of variations with a rich history of mathematical rivalry associated with it, the geodesic problem would go on to outpace it in terms of importance with the development of differential geometry. Eventually geodesics would become an essential part in our understanding of motion under gravity with the advent of general relativity \cite{weinberg1972gravitation}. Here we look at the simple problem of finding the shortest path on unit a sphere given two points $(1, \theta_0, \phi_0), (1, \theta_1, \phi_1)$ (in spherical polar coordinates) on it by setting,
\begin{equation}
\begin{aligned}
    &\Omega=[\theta_0, \theta_1],\;X=W^{1, 2}(\Omega; [0, 2\pi)),\; W=\mathbb R\\
    &f(u) = \int_{\theta_0}^{\theta_1}\sqrt{1+\left(\sin\theta \frac{du}{d\theta}\right)^2}\;\,d\theta\\
    & g(u) = \sqrt{(u(\theta_0)-\phi_0)^2 + (u(\theta_1)-\phi_1)^2}\label{eq:gs}
\end{aligned}
\end{equation}
If $u^*$ is the solution then $(1,\theta, u^*(\theta))$ gives us a parametrization for the geodesic curve.
\subsection{Grad-Shafranov equation} Grad-Shafranov equation is an elliptic partial differential equation describing the poloidal flux under ideal magnetohydrodynamics for a 2D plasma \cite{smithaaxisymmetric}.  Modelling the plasma equilibrium is an important aspect of designing magnetic confinement devices like tokamaks in the field of nuclear fusion. Although originally used for axis-symmetric tokamaks, the Grad-Shafranov equation has been analyzed for non-axis symmetric magnetohydrodynamic equilibrium as well \cite{burby2020generalized}. In 1968 Solov'ev derived a family of analytic solutions for the Grad-Shafranov equation under the assumption that there is distributed toroidal
current filling all space \cite{xu2019vacuum} and since then these Solev'ev solutions have become an import benchmarking tool for plasma equilibrium codes \cite{johnson1979numerical}. Below we describe the Grad-Shafranov equation, this specific version can also be found in \cite{xu2019vacuum}. 
\begin{equation}
\begin{aligned}
    &\frac{\partial^2 u}{\partial z^2}+r\frac{\partial}{\partial r}\left(\frac{1}{r}\frac{\partial u}{\partial r}\right) =ar^2 + bR^2,\qquad (r,z)\in\Omega= [0.9R, 1.1R]\times[-0.1R, 0.1R]\\
    & u(r, z) = \frac{1}{2}(b+c_0)R^2z^2+c_0R\zeta z^2+\frac{1}{2}(a-c_0)R^2\zeta^2, \qquad (r, z)\in{\partial\Omega}\\
    &\text{\rm where }\zeta =\frac{r^2-R^2}{2R}, R=1.0,\, a = 1.2,\,b=-1.0,\, c_0=1.1
\end{aligned}
\end{equation}
In order to cast this problem into the format of \eqref{eq:con-opt}, we set
\begin{equation}
\begin{aligned}
    &X = W^{1,2}(\Omega;\mathbb R),\;W=L^{2}(\partial\Omega; \mathbb R)\\
    &f(u) = \int_{-0.1R}^{0.1R}\int_{0.9R}^{1.1R}\left(\frac{\partial^2 u}{\partial z^2}+r\frac{\partial}{\partial r}\left(\frac{1}{r}\frac{\partial u}{\partial r}\right) -ar^2 - bR^2\right)^2\,dr\,dz\\
    &g(u): (r,z)\mapsto \frac{1}{2}(b+c_0)R^2z^2+c_0R\zeta z^2+\frac{1}{2}(a-c_0)R^2\zeta^2 - u(r, z)
\end{aligned}\label{eq:GS}
\end{equation}
\subsection{Beltrami fields} Beltrami fields are special vector fields that are eigenfunctions of the curl operator. They play an important role in fluid dynamics as steady solutions to the Euler equation \cite{aris2012vectors}. In this problem we ask, given Beltrami boundary data, what is the magnetic field of least energy in a 3D volume? Gauss's law \cite{jackson1999classical} dictates that we have to take the nondivergence of magnetic fields into account which can be done in multiple ways while formulating our question, either as a part of the Hilbert space $X$ (since divergence is a linear operator) or as an addition to the boundary condition $g$. Here we choose to impose Gauss's law as a part of the Hilbert space $X$.
\begin{equation}
    \begin{aligned}
&\Omega = \left[-\frac{1}{2}, \frac{1}{2}\right]^3,\; X=\overline{\{u\in W^{1,2}(\Omega; \mathbb R^3):\nabla\cdot u=0\}},\;W=L^{2}(\partial\Omega;\mathbb R^3)\\
    &f(u) = \frac{1}{2}\int_{-\frac{1}{2}}^{\frac{1}{2}}\int_{-\frac{1}{2}}^{\frac{1}{2}}\int_{-\frac{1}{2}}^{\frac{1}{2}}|u(x, y, z)|^2\,dx\,dy\,dz\\
    & g(u): (x, y, z)\mapsto u(x, y, z) - \begin{bmatrix}
\sin(z) + \cos(y) \\
\sin(x) + \cos(z) \\
\sin(y) + \cos(x) \\
\end{bmatrix}
    \end{aligned}
\end{equation}
Unlike the other problems stated here, this problem is \textit{manufactured} and has no direct practical applications but nevertheless serves as an interesting toy problem.

%% file: var-tex/method.tex
Before discussing our algorithm for solving the problems stated in section~\ref{sec-prob}, we briefly look at constrained optimization algorithms for finite dimensional problems and their infinite dimensional analogues since they illustrate the guiding principles that will help us devise our own algorithm.

\subsection{Constrained optimization algorithms in finite dimensions}
Unconstrained optimization problems are typically easier to solve than constrained optimization problems and they are often solved using variants of gradient descent or Newton's method \cite{bonnans2006numerical}, \cite{jorge2006numerical}. Therefore, in order to solve constrained problems we often transform them into unconstrained problems first. When $X, W$ are finite dimensional, in order to solve problem~\eqref{eq:con-opt} we can convert it from a constrained optimization problem to a sequence of unconstrained subproblems as follows,
\begin{align}
    u_k=\underset{u\in X}\arginf \mathcal L(u, \mu_k)\stackrel{\rm def}{=}f(u) + \frac{\mu_k}{2}|g(u)|_W^2,\;k=1,2,\cdots\label{eq:seq-uncon}
\end{align}
where $|\cdot|_W$ and  denotes the canonical norm on $W$ and $\{\mu_k\}_{k=1}^\infty$ is a positive increasing sequence such that $\mu_k\uparrow\infty$ and $u_k$ is the exact global solution to the $k$-th subproblem. It can be shown that every limit point of $\{u_k\}_{k=1}^\infty$ is a solution to the original constrained problem, for a proof see theorem~17.1 in \cite{jorge2006numerical} or for a local version of the statement see theorem~1 in \cite{polyak1971convergence}. The strategy of using subproblems \eqref{eq:seq-uncon} to solve \eqref{eq:con-opt} is known as the quadratic \textit{penalty method}. In case we have access to only approximate solutions to the subproblems then limit points of these approximate solutions might be infeasible or they might only satisfy the first order KKT condition  \cite{kuhn2013nonlinear}, \cite{gordon2012karush}, \cite{boltyanski1998geometric} rather than being global minimizers. Moreover, the Hessian of the unconstrained objective function $\mathcal L$ becomes ill-conditioned as $\mu_k\uparrow\infty$. If we attempt to find an approximate solution to \eqref{eq:seq-uncon} by trying to satisfy the first order condition using Newton's method, we quickly run into significant numerical errors when $\mu_k$ is large. For an excellent discussion of the nuances associated with the penalty method see chapter~17 in \cite{jorge2006numerical}. Typically the convergence rate of the quadratic penalty method is $O(k^{-\frac{1}{2}})$ but for strongly convex problems it increases to $O(k^{-1}) $  \cite{li2017convergence}, \cite{polyak1971convergence}.

In order to avoid ill-conditioning we can modify our subproblems as follows,
\begin{align}
    u_k=\underset{u\in X}\arginf\mathcal L_A(u, \mu_k, \lambda_k)\stackrel{\rm def}{=} f(u) + \frac{\mu_k}{2}|g(u)|_W^2 + \langle\lambda_k, g(u)\rangle_W,\;k=1,2,\cdots\label{eq:seq-uncon-al}
\end{align}
where $\langle\cdot\rangle_W$ is the canonical inner product on $W$, $\{\mu_k\}_{k=1}^\infty$ is a positive, nondecreasing sequence but not necessarily unbounded and $\lambda_k$ follows the update rule,
\begin{align}
    \lambda_{k+1} = \lambda_k + \mu_k g(u_k)\label{eq:mul-update}
\end{align}
This update rule is the consequence of an attempt to satisfy the first order condition for optimality. In this setting, $\mu_k$ and $\lambda_k$ play the roles of the penalty factor for deviating from the constraint and the Lagrange multiplier respectively. This method is known as the \textit{augmented Lagrangian method}. It can be shown that, under suitable conditions, if $\lambda_k$ converges to $\lambda^*$ then $\exists \,\mu^*>0$ such that for $\mu\ge\mu^*$, any local solution to the original constrained problem is a local minimizer of $\mathcal L_A(\cdot, \mu, \lambda^*)$, see theorem~17.5 in \cite{jorge2006numerical} or for a global version of this statement see theorem~5.2 in \cite{birgin2014practical}. This gives the augmented Lagrangian method a strong theoretical foundation but in practice we might only have approximate knowledge of $\lambda^*$. In such a case i.e. when $\lambda_k$ is close to $\lambda^*$, it can be shown that a local minimizer of $\mathcal L_A(\cdot, \mu, \lambda_k)$ solves the original constrained problem for large enough $\mu$, see theorem~17.6 in \cite{jorge2006numerical} or proposition~4.2.3 in \cite{bertsekas1995athena}. These results show that the augmented Lagrangian method can approximately solve \eqref{eq:con-opt} when either the penalty $\mu$ is large or we have good knowledge of the optimal Lagrange multiplier $\lambda^*$. The appeal of the augmented Lagrangian method therefore lies in the possibility that we can replace the requirement that $\mu_k\uparrow\infty$ with the convergence of the Lagrange multiplier $\lambda_k$ thus avoiding the ill-conditioning of the Hessian and all the numerical difficulties that arise because of it.

\subsection{Constrained optimization algorithms in infinite dimensions} When $X$ is infinite dimensional and $f(\cdot),\langle g(\cdot), g(\cdot)\rangle_W$ are lower-semicontinuous functionals, limit points of the exact global solutions to the subproblems are solutions to the original constrained problem, for a proof see theorem~1 in \cite{dussault2007penalty} and for a treatment of the penalty method on general topological spaces see \cite{fiacco1969generalized}. The augmented Lagrangian method has also been extended in many different scenarios where $X$ is an infinite dimensional Hilbert space by Ito and Kunisch \cite{ito1990augmented}, \cite{ito1990augmentedvar}, \cite{ito2008lagrange}. More recently the case when $X$ is an infinite dimensional Banach space has been considered by Kanzow et al \cite{kanzow2018augmented}. If we assume that problem~\eqref{eq:con-opt} has a solution, $f, g$ are twice continuously Fréchet differentiable near the solution, derivative of $g$ at the solution is surjective, a Lagrange multiplier exists for this solution, $f$ is weakly lower-semicontinuous and $g$ maps weakly convergent sequences to weakly convergent sequences then we can prove that the augmented Lagrangian subproblems have local solutions, the Lagrange multiplier $\lambda_k$ converges to $\lambda^*$ and these local solutions converge to a local solution of \eqref{eq:con-opt}. For an in-depth look at the technical details, we refer the reader to chapter~3 of \cite{ito2008lagrange}. 

\subsection{Deep learning variants for infinite dimensional algorithms}
The first challenge in implementing these algorithms is representing elements of $X$ and $W$ when they are infinite dimensional. A direct approach to do this would be to represent an element of $X$ as a neural network $u_{\eta}^\mathcal{A}$  and in case $W$ is infinite dimensional, we can represent the Lagrange multiplier as another network $\lambda_\xi^\mathcal{B}$ where $\eta, \xi$ represent the trainable parameters of the networks and $\mathcal A, \mathcal B$ represent the structure or architecture of the networks. Universal approximation theorems~\cite{pinkus1999approximation}, \cite{de2021approximation} imply that with appropriately chosen $\mathcal A, \mathcal B$ we might be able to sufficiently approximate the solutions to the subproblems. Suppose the dimensions of $\eta, \xi$ or the number of trainable parameters are $a, b$ respectively. Then the subproblems in the penalty algorithm can be rewritten as,
\begin{align}
    \eta_k=\underset{\eta\in \mathbb R^a}\arginf \mathcal L(u_{\eta}^{\mathcal A}, \mu_k)=f(u_{\eta}^{\mathcal A}) + \frac{\mu_k}{2}|g(u_{\eta}^{\mathcal A})|_W^2,\;k=1,2,\cdots\label{eq:seq-uncon-net}
\end{align}
Similarly, the subproblems in the augmented Lagrangian algorithm can be rewritten as, 
\begin{align}
    \eta_k=\underset{\eta\in \mathbb R^a}\arginf \mathcal L_A(u_{\eta}^{\mathcal A}, \mu_k, \lambda_{\xi_k}^{\mathcal B})=f(u_{\eta}^{\mathcal A}) + \frac{\mu_k}{2}|g(u_{\eta}^{\mathcal A})|_W^2 + \langle\lambda_{\xi_k}^{\mathcal B}, g(u_{\eta}^{\mathcal A})\rangle_W,\;k=1,2,\cdots\label{eq:seq-uncon-al-net}
    \end{align}
When $W$ is infinite dimensional, the Lagrange multiplier update rule can be rewritten as,
\begin{align}
\lambda_{\xi_{k+1}}^{\mathcal B} = \lambda_{\xi_k}^{\mathcal B} + \mu_k g(u_{\eta_k}^{\mathcal A})\label{eq:mul-update-net}
\end{align}
 Note that, with this rewriting our infinite-dimensional subproblems have become finite dimensional since $a$ is finite. The update rule \eqref{eq:mul-update-net} can not be implemented directly since the Lagrange multipliers are functions rather than finite dimensional vectors in this scenario. Therefore, we try to find the optimal $\xi_{k+1}$ that makes the left hand side of \eqref{eq:mul-update-net} functionally or in an $L^2$ sense, equal to the right hand side by solving the following optimization problem, 
\begin{align}
    \xi_{k+1}=\underset{\xi\in\mathbb R^b}\arginf\left|\lambda_{\xi}^{\mathcal B}-\lambda_{\xi_{k}}^{\mathcal B}-\mu_k g(u_{\eta_k}^{\mathcal A})\right|_W^2\label{eq:xi-update}
\end{align}
If we solve $K$ subproblems then we approximate our final solution as $u_{\eta_K}$. In order to solve each subproblem in \eqref{eq:seq-uncon-net} and \eqref{eq:seq-uncon-al-net} we use our solution to the last subproblem as an initial guess and perform gradient descent. In order to solve \eqref{eq:xi-update} we use $\xi_k$ as an initial guess since we expect the sequence $\lambda_{\xi_k}^{\mathcal B}$ to converge. The selection of $\mu_k$ is an important part of the algorithm but no general purpose techniques for this selection are available in the literature. Larger $\mu_k$ results in better theoretical convergence rates while deteriorating the numerical estimates at same time, see section~3.3 in \cite{ito2008lagrange} for comments on this topic. The update rule in \eqref{eq:mul-update} can be modified in different ways to achieve better estimates when one has some extra information about the problem \eqref{eq:con-opt}, see the ALM algorithm in \cite{ito2008lagrange} for example. But in practice such information is nearly impossible to come by and therefore we will stick to the simple update rule in \eqref{eq:mul-update}. The deep learning variants of the penalty method, augmented Lagrangian method when $W$ is finite dimensional and augmented Lagrangian method when $W$ is infinite dimensional can found in algorithms \ref{algo:dl-penalty} (${\rm P}^\infty$), \ref{algo:dl-al-finite} (${\rm AL}_{\rm F}^\infty$), \ref{algo:dl-al-infinite} (${\rm AL}_\infty^\infty$) respectively. 

\begin{algorithm}[!ht]
\caption{${\rm P}^\infty$: Infinite dimensional penalty algorithm}
\label{algo:dl-penalty}
\begin{algorithmic}[1]
    \STATE Choose architecture $\mathcal A$, penalty factor sequence $\{\mu_k\}_{k=1}^\infty$, adaptive learning rate $\{\delta_{k, j}\}_{k, j=1}^\infty$ and stopping criteria $\{P_k\}_{k=1}^\infty, \{Q_{k, j}\}_{k, j=1}^\infty$
    \STATE $k\leftarrow0$
    \WHILE{stopping criterion $P_k$ is not met}
        \STATE $k \leftarrow k+1$
        \IF{$k=1$}
                \STATE Initialize $\eta$ randomly
            \ELSE
                \STATE Initialize $\eta\leftarrow \eta_{k-1}$
        \ENDIF
        \STATE $j\leftarrow 1$
        \WHILE{stopping criterion $Q_{k, j}$ is not met}
            \STATE $L\leftarrow f(u_{\eta}^{\mathcal A})+\frac{\mu_k}{2}|g(u_{\eta}^{\mathcal A})|^2_W$
            \STATE $\eta\leftarrow\eta-\delta_{k, j}\nabla_\eta L$
            \STATE $j\leftarrow j+1$
        \ENDWHILE
        \STATE $\eta_k\leftarrow\eta$
    \ENDWHILE
    \STATE $u_{\eta_k}^{\mathcal A}$ is our approximate solution to \eqref{eq:con-opt}
\end{algorithmic}
\end{algorithm}

\begin{algorithm}[!ht]
\caption{${\rm AL}^\infty_{\rm F}$: Infinite dimensional augmented Lagrangian algorithm when $W$ is finite dimensional}
\label{algo:dl-al-finite}
\begin{algorithmic}[1]
    \STATE Choose architectures $\mathcal A$, penalty factor sequence $\{\mu_k\}_{k=1}^\infty$, adaptive learning rate $\{\delta_{k, j}\}_{k, j=1}^\infty$ and stopping criteria $\{P_k\}_{k=1}^\infty, \{Q_{k, j}\}_{k, j=1}^\infty$
    \STATE $k\leftarrow0$
    \WHILE{stopping criterion $P_k$ is not met}
        \STATE $k \leftarrow k+1$
        \IF{$k=1$}
                \STATE Initialize $\eta$ randomly
        \ELSE
                \STATE Initialize $\eta\leftarrow \eta_{k-1}$
        \ENDIF
        \STATE $j\leftarrow 1$
        \WHILE{stopping criterion $Q_{k, j}$ is not met}
            \STATE $L\leftarrow f(u_{\eta}^{\mathcal A})+\frac{\mu_k}{2}|g(u_{\eta}^{\mathcal A})|^2_W + \langle\lambda_{\xi_k}^{\mathcal B}, g(u_{\eta}^{\mathcal A})\rangle_W$
            \STATE $\eta\leftarrow\eta-\delta_{k, j}^A\nabla_\eta L$
            \STATE $j\leftarrow j+1$
        \ENDWHILE
        \STATE $\eta_k\leftarrow\eta$
        \STATE $\xi_{k+1}\leftarrow\xi_k + g(u_{\eta_k}^{\mathcal A})$
    \ENDWHILE
    \STATE $u_{\eta_k}^{\mathcal A}$ is our approximate solution to \eqref{eq:con-opt}
\end{algorithmic}
\end{algorithm}

\begin{algorithm}[!ht]
\caption{${\rm AL}_\infty^\infty$: Infinite dimensional augmented Lagrangian algorithm when $W$ is infinite dimensional}
\label{algo:dl-al-infinite}
\begin{algorithmic}[1]
    \STATE Choose architectures $\mathcal A$, $\mathcal B$, penalty factor sequence $\{\mu_k\}_{k=1}^\infty$, adaptive learning rates $\{\delta_{k,j}^A\}_{k,j=1}^\infty$, $\{\delta_{k, j}^B\}_{k,j=1}^\infty$ and stopping criteria $\{P_k\}_{k=1}^\infty, \{Q_{k,j}^A\}_{k, j=1}^\infty, \{Q_{k,j}^B\}_{k, j=1}^\infty$
    \STATE $k\leftarrow0$
    \WHILE{stopping criterion $P_k$ is not met}
        \STATE $k \leftarrow k+1$
        \IF{$k=1$}
                \STATE Initialize $\eta$ randomly
                \STATE Initialize $\xi$ randomly
            \ELSE
                \STATE Initialize $\eta\leftarrow \eta_{k-1}$
        \ENDIF
        \STATE $j\leftarrow 1$
        \WHILE{stopping criterion $Q_{k, j}^A$ is not met}
            \STATE $L\leftarrow f(u_{\eta}^{\mathcal A})+\frac{\mu_k}{2}|g(u_{\eta}^{\mathcal A})|^2_W + \langle\lambda_{\xi_k}^{\mathcal B}, g(u_{\eta}^{\mathcal A})\rangle_W$
            \STATE $\eta\leftarrow\eta-\delta_{k,j}^A\nabla_\eta L$
            \STATE $j\leftarrow j+1$
        \ENDWHILE
        \STATE $\eta_k\leftarrow\eta$
        \STATE $j\leftarrow 1$
        \WHILE{stopping criterion $Q_{k, j}^B$ is not met}
            \STATE $L_\lambda\leftarrow \left|\lambda_{\xi}^{\mathcal B}-\lambda_{\xi_{k}}^{\mathcal B}-\mu_{k} g(u_{\eta_{k}}^{\mathcal A})\right|_W^2$
            \STATE $\xi\leftarrow\xi-\delta_{k,j}^B\nabla_\xi L_\lambda$
            \STATE $j\leftarrow j+1$
        \ENDWHILE
        \STATE $\xi_{k+1}\leftarrow\xi$
    \ENDWHILE
    \STATE $u_{\eta_k}^{\mathcal A}$ is our approximate solution to \eqref{eq:con-opt}
\end{algorithmic}
\end{algorithm}

%% file: var-tex/results.tex
In this section we describe the results along with the specific details of the algorithms used for each problem. We solve every problem with both penalty and augmented Lagrangian methods.
\begin{itemize}
    \item \textbf{Architecture}: We use two different types of architecture in our experiments, see appendix~\ref{ssec-arch} for a description of these types and refer to table~\ref{tab:network} for details of $\mathcal A, \mathcal B$ used in the experiments. We use the same architecture to represent the approximate solutions of penalty and augmented Lagrangian algorithms.
    \item \textbf{Stopping criteria}: Although sophisticated stopping criteria such as the norm of the gradient of the objective function in the subproblem falling below a pre-selected threshold, can be used for algorithms~\ref{algo:dl-penalty}, \ref{algo:dl-al-finite}, \ref{algo:dl-al-infinite}, here we stop the loops after a pre-selected number of iterations is reached.
    \item \textbf{Number of gradient descent steps}: We denote this number of iterations with $P$ for the outer loops in algorithms~\ref{algo:dl-penalty}, \ref{algo:dl-al-finite}, \ref{algo:dl-al-infinite}, $Q$ for the inner loops in algorithms~\ref{algo:dl-penalty}, \ref{algo:dl-al-finite} and $Q^A,Q^B$ for the first and second inner loops in the algorithm~\ref{algo:dl-al-infinite} respectively. We define $E$ to be the total number of gradient descent steps used. Therefore, $E=PQ$ for algorithms~\ref{algo:dl-penalty} and \ref{algo:dl-al-finite} and $E=P(Q^A+Q^B)$ for algorithm~\ref{algo:dl-al-infinite}. For each problem we use the same total gradient descent steps $E$ for both the penalty and the augmented Lagrangian algorithms. To facilitate this, when $W$ is finite dimensional we use the same $Q$ for algorithms~\ref{algo:dl-penalty} and \ref{algo:dl-al-finite} and when $W$ is infinite dimensional we set $Q^A=Q^B=Q/2$. We use the popular Adam optimizer \cite{kingma2014adam} to perform the gradient descent steps.
    \item\textbf{Learning rate}: We use an initially oscillating and finally decaying learning rate $\delta$ that depends on 7 distinct hyperparameters. The oscillatory nature of $\delta$ as seen in figure~\ref{fig:learning-rate}, is employed to essentially rejuvenate the previously decaying learning rate every time we start  an inner loop in the algorithms. For details of this learning rate $\delta$, see appendix~\ref{ssec-rate}. 
\begin{figure}[!ht]
    \centering
\includegraphics[scale=0.5
]{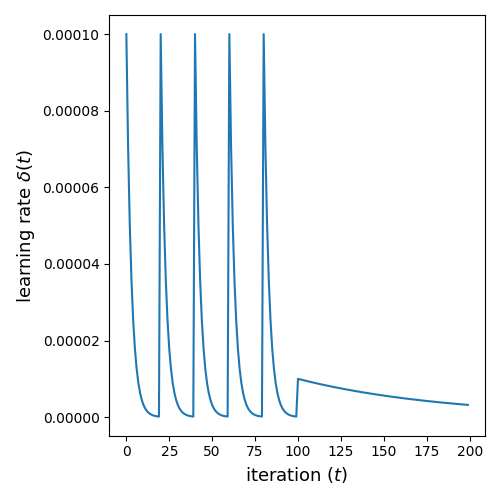}
    \caption{Example behavior of the oscillating learning rate $\delta$}
    \label{fig:learning-rate}
\end{figure}
While using algorithms~\ref{algo:dl-penalty} and \ref{algo:dl-al-finite} we set,
\begin{align}
    \delta_{k, j}=\delta((k-1)Q+j)
\end{align}
and while using algorithm~\ref{algo:dl-al-infinite} we use,
\begin{align}
    \delta_{k, j}^A = \delta((k-1)Q^A+j)\\
    \delta_{k, j}^B = \delta((k-1)Q^B+j)
\end{align}
\item\textbf{The penalty factors}: We use a stopped geometric sequence as our $\mu_k$,
\begin{align}
    \mu_k=\min\{\mu_1r^{k-1},\mu_{\max}\}\label{eq:mu}
\end{align}
The exact values of $\mu_1,\mu_{\max}, r$ for various problems can be found in table~\ref{tab:mu}.
\item \textbf{Computation of functionals}: To compute the functional $f$ and when $W$ is infinite dimensional, the functional $|\cdot|_W$, we use either Gauss-Legendre quadrature (in 1 or 2 dimensions) or a Monte Carlo estimate (in 3 dimensions).
\item\textbf{Errors}: We evaluate our algorithms using three different kinds of errors produced.
If $\hat u$ is the solution produced by an algorithm and $u^{\rm true}$ is the true solution of the problem, we define the absolute error to be an weighted $L^2$-norm  of $\hat u-u^{\rm true}$,
\begin{align}
{\rm absolute\;error}=\sqrt{\frac{\int_{\Omega}(\hat u-u^{\rm true})^2\,dV}{\int_{\Omega}\,dV}}\label{eq:abs-err}
\end{align}
where $dV$ denotes a volume element in $\Omega$. We define the relative objective error to be the relative error in the value of the objective function, 
\begin{align}
{\rm relative\; objective\; error}=\left|\frac{f(\hat u)-f(u^{\rm true})}{f(u^{\rm true})}\right|\label{eq:rel-err}
\end{align}
\end{itemize}
Lastly, we define the constraint error to be how closely $\hat u$ satisfies the constraint,
\begin{align}
    {\rm constraint\; error}= \frac{|g(\hat u)|_W}{Z}\label{eq:cons-err}
\end{align}
We set $Z=1$ for the geodesic problem~\eqref{eq:gs} and $Z=\sqrt{|\partial\Omega|}$ for all the other problems. The use of normalization constants in \eqref{eq:abs-err}, \eqref{eq:cons-err} results in errors that are akin to RMSE. With these definitions we are now ready to present the results.     

\subsection{The minimal surface problem}\label{ssec:res-ms}
\begin{figure}[!ht]
    \centering
\includegraphics[scale=0.4]{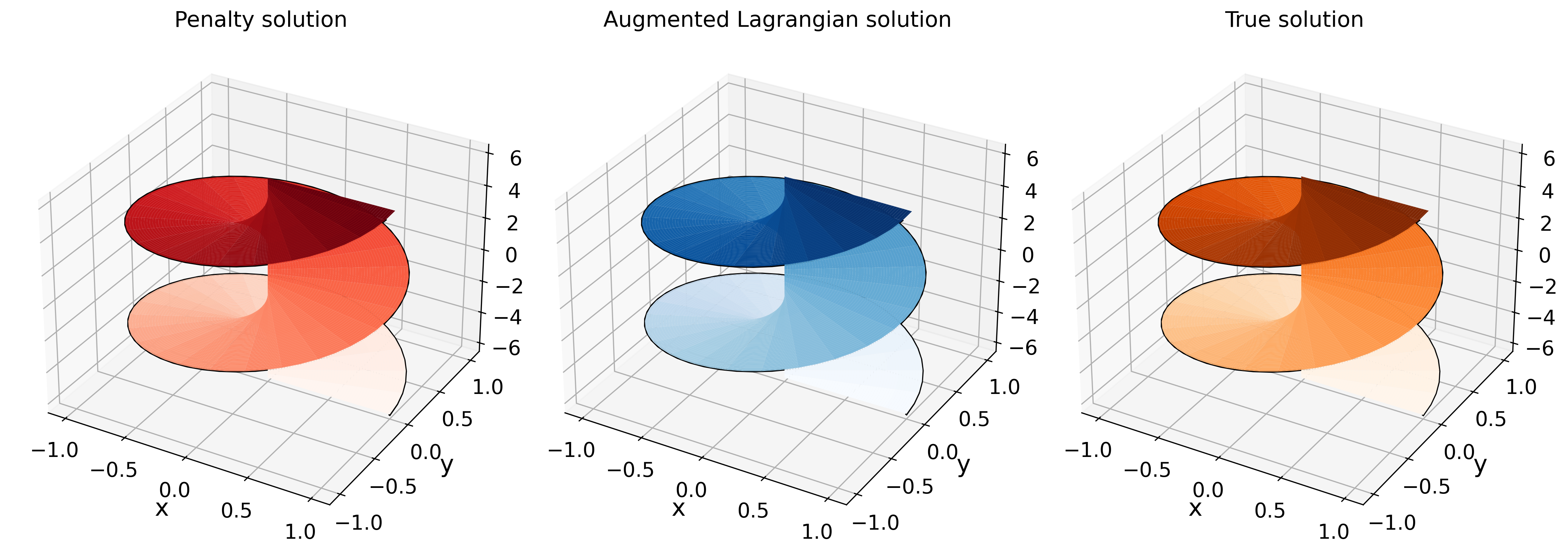}
    \caption{Solutions to the minimal surface problem. Darker color implies higher $u$ value.}
    \label{fig:helicoid-surface}
\end{figure}

\begin{figure}[!ht]
    \centering
\includegraphics[scale=0.32]{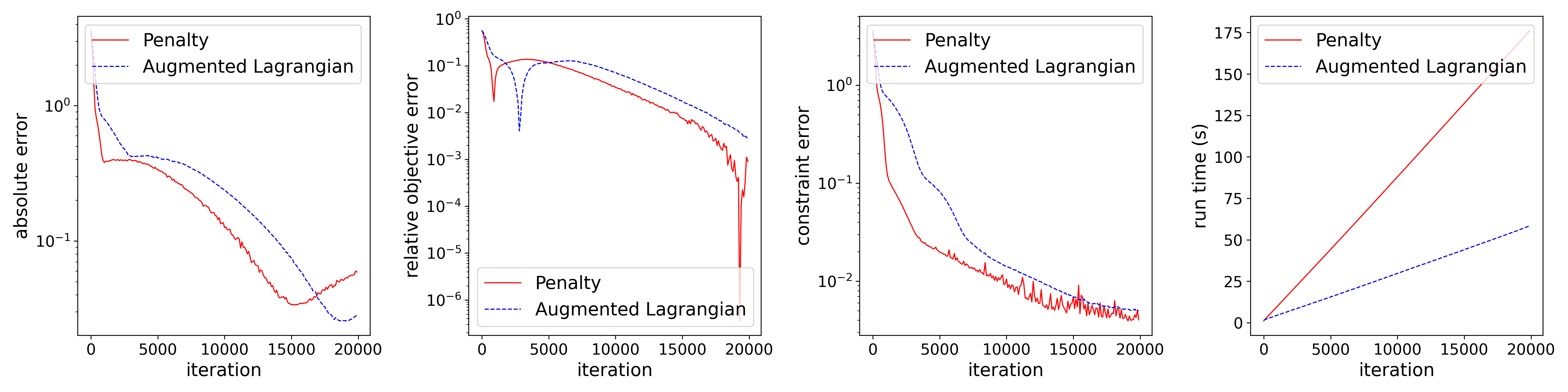}
    \caption{Errors and run times for the minimal surface problem as functions of gradient descent steps. The errors have been plotted in a semilog fashion. All quantities have been plotted every 100 steps.}
    \label{fig:helicoid-error}
\end{figure}

Figure~\ref{fig:helicoid-surface} shows the true and approximate solutions to the minimal surface problem in Cartesian coordinates. The true solution for this problem is a helicoid. 
We use $E=20000$ total gradient steps and $P=1000$ subproblems for this problem. This setup implies to calculate solutions to the subproblems, we use $Q=20$ gradient descent steps for the penalty algorithm (${\rm P}^\infty$) and $Q^A=Q^B=10$ gradient steps for the augmented Lagrangian algorithm (${\rm AL}^\infty_\infty$). Figure~\ref{fig:helicoid-error} shows various errors as functions of gradient descent steps for this problem. In terms of algorithms~\ref{algo:dl-penalty}, \ref{algo:dl-al-finite}, \ref{algo:dl-al-infinite}, iteration in figure~\ref{fig:helicoid-error} can be understood as $(k-1)E/P+j$. The errors and run times have been plotted every 100 steps.
Both methods are able to produce good approximations but the penalty method fluctuates more during training compared to its counterpart for this problem. Looking at table~\ref{tab:network} we see that sizes of the networks representing the solution and the Lagrange multiplier $a, b$ are close to each other. But solving the problem \eqref{eq:xi-update} is computationally much cheaper than solving the problem \eqref{eq:seq-uncon-al-net}. This is a typical scenario and is reflected in the run times of the algorithms in figure~\ref{fig:helicoid-error} where the augmented Lagrangian is $3.016$ times faster than its counterpart.

\subsection{Geodesics on a surface}
\begin{figure}[!ht]
    \centering
\includegraphics[scale=0.4]{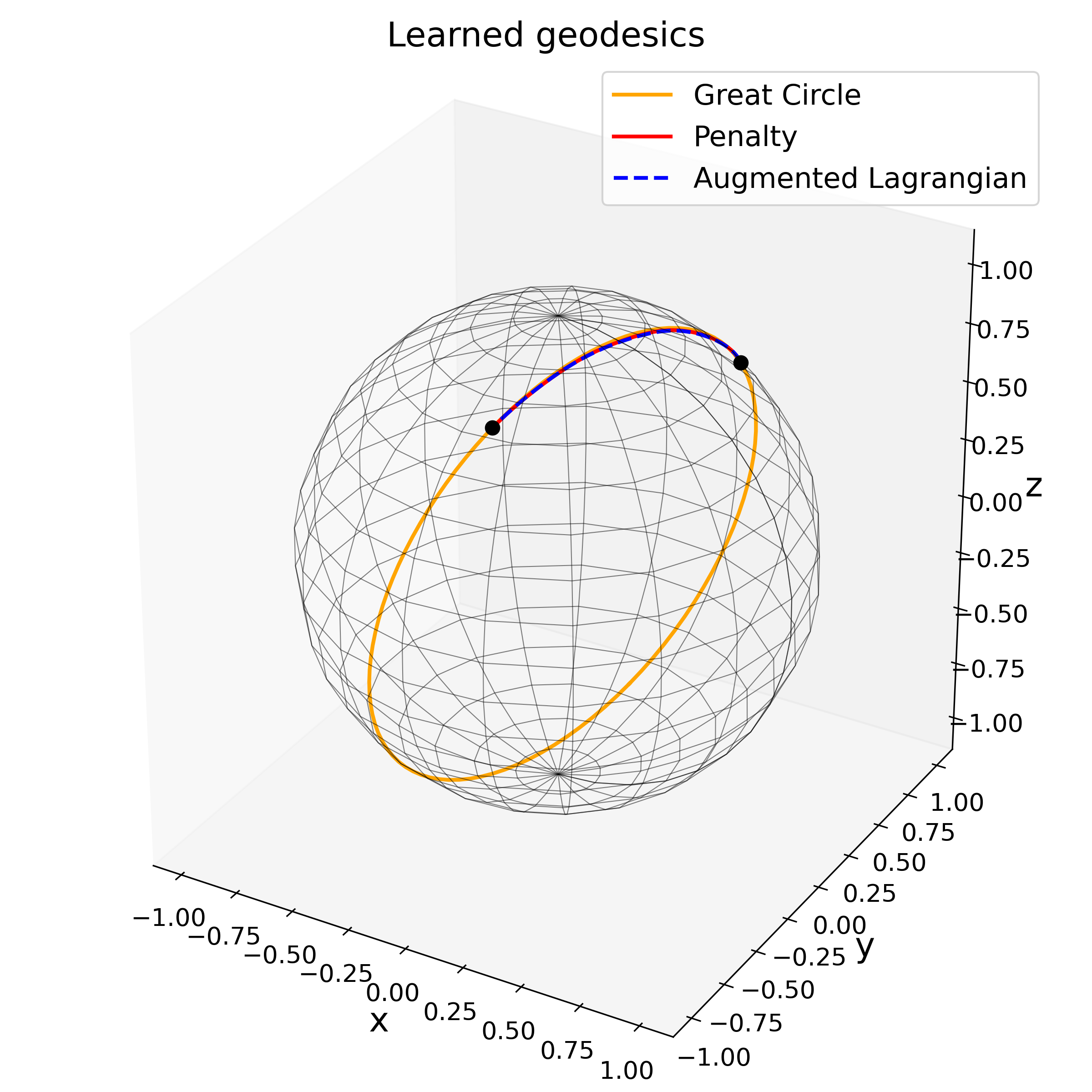}
    \caption{Solutions to the geodesic problem. The distance between the black dots is being minimized.}
    \label{fig:geo}
\end{figure}

\begin{figure}[!ht]
    \centering
\includegraphics[scale=0.32]{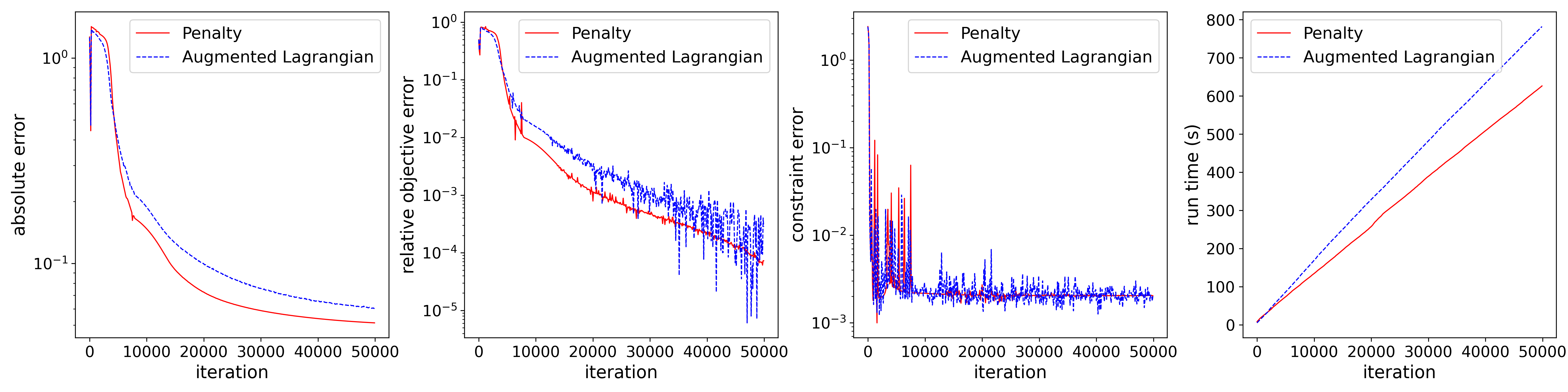}
    \caption{Errors and run times for the geodesic problem as functions of gradient descent steps. The errors have been plotted in a semilog fashion. All quantities have been plotted every 100 steps.}
    \label{fig:geo-error}
\end{figure}

Figure~\ref{fig:geo} shows the approximate solutions to the geodesic problem. The true solution to this problem is the arc between the given points that lies on the great circle connecting them.
We use $E=50000$ and $P=2500$ for this problem which implies the number of gradient steps used to solve subproblems for the both the penalty (${\rm P}^\infty$) and the augmented Lagrangian (${\rm AL}^\infty_{\rm F}$) algorithms is $Q=20$. Figure~\ref{fig:geo-error} shows the errors and run times for this problem as a function of gradient steps.
 Although the absolute error decreases without significant fluctuations, ${\rm AL}^\infty_{\rm F}$ shows considerably more fluctuations in other errors compared to its counterpart during training. Since \eqref{eq:seq-uncon-al-net} is computationally more expensive to solve than \eqref{eq:seq-uncon-net} and we use the same number of gradient descent steps to solve both of them, the augmented Lagrangian algorithm is slower in this case as seen in figure~\ref{fig:geo-error}.

\begin{figure}[!ht]
    \centering
\includegraphics[scale=0.4]{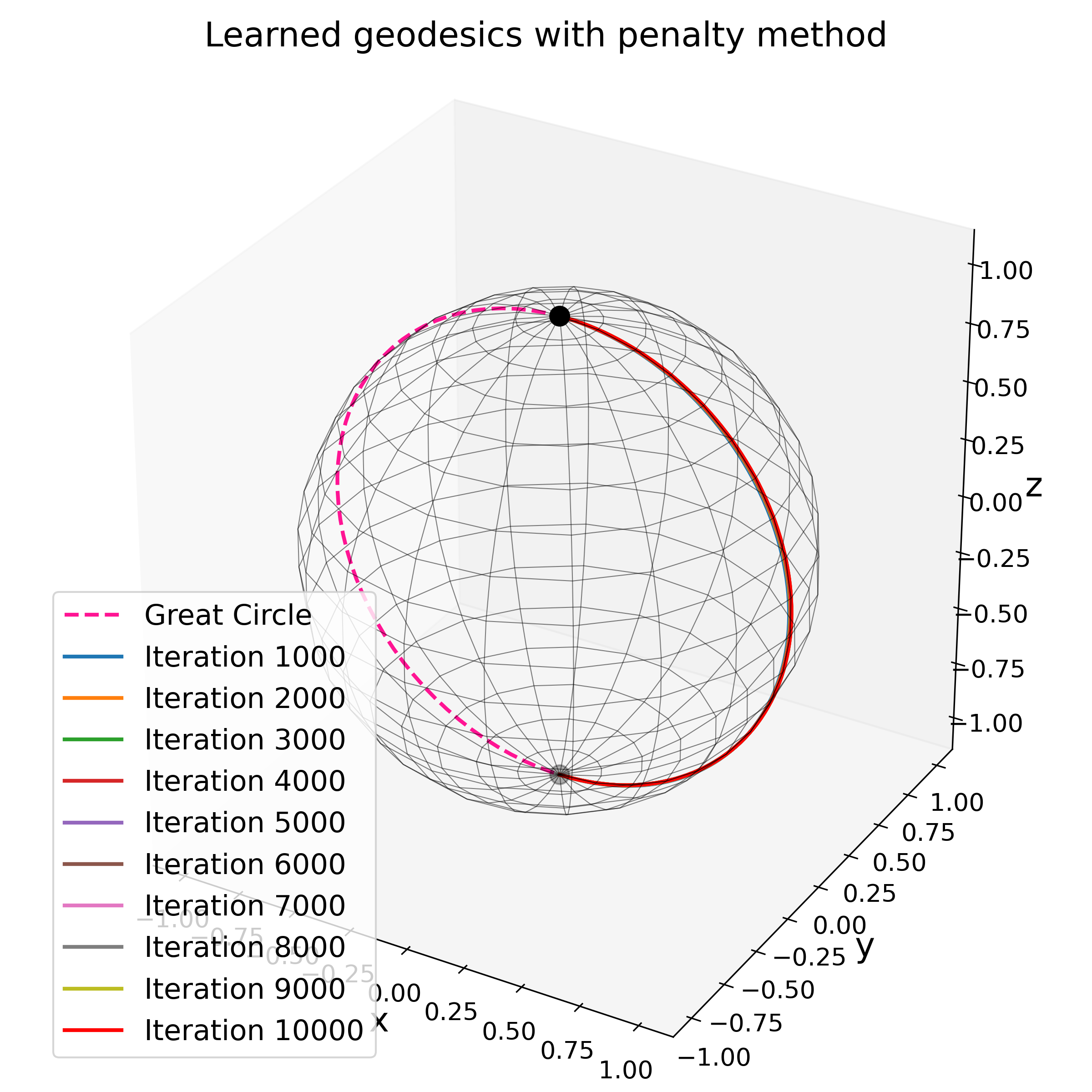}
    \caption{Solutions to the geodesic problem when the points (black dots) are antipodal}
    \label{fig:geo-ap}
\end{figure}

In figure~\ref{fig:geo} the points on the sphere are chosen to be nonantipodal, leading to a unique solution to the geodesic problem. In case these points are antipodal there are infinitely many great circles that connect them, leading to infinitely many solutions. The solution set in this degenerate case is homeomorphic to a connected 1D manifold. Figure~\ref{fig:geo-ap} shows the geodesics learned with the penalty method for such a degenerate problem. As can be seen, the solution produced by the penalty algorithm (${\rm P}^\infty$) quite surprisingly does not hop from one great circle to another but rather stays on a single great circle with more and more gradient descent steps. Even though the true minima of the original problem lie on a connected manifold, the discretized version of the problem where we look for $\eta\in\mathbb R^a$, might have a different distribution of minima, where possibly each minimum can be separated from any other minimum by open sets. The geometry of the loss landscape in the subproblems is an interesting topic and requires further investigation.
\subsection{Grad-Shafranov equation}
\begin{figure}[!ht]
    \centering
\includegraphics[scale=0.32]{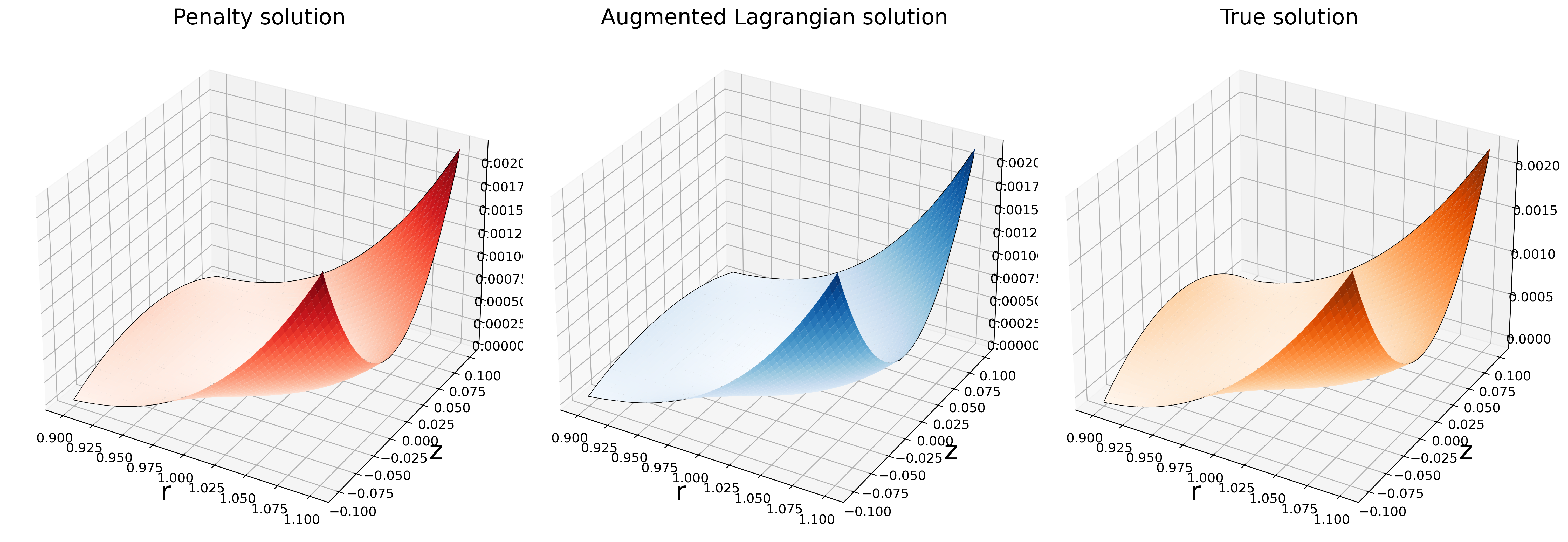}
    \caption{Solutions to the Grad-Shafranov equation. Darker color implies higher $u$ value.}
    \label{fig:gs-surface}
\end{figure}

\begin{figure}[!ht]
    \centering
\includegraphics[scale=0.32]{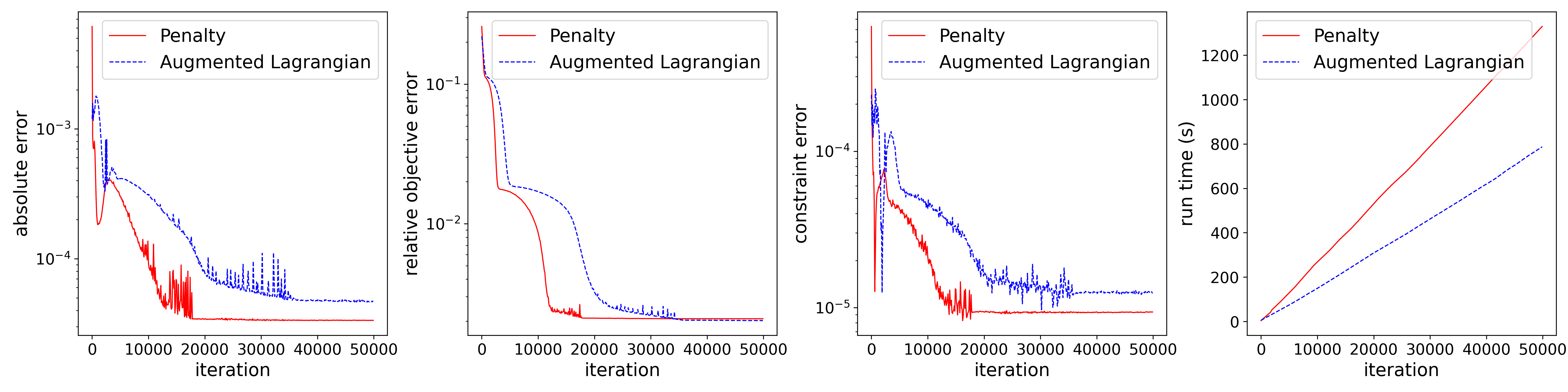}
    \caption{Errors and run times for the Grad-Shafranov problem as functions of gradient descent steps. The errors have been plotted in a semilog fashion. All quantities have been plotted every 100 steps.}
    \label{fig:gs-error}
\end{figure}

Figure~\ref{fig:gs-surface} shows the approximate and true solutions for the Grad-Shafranov equation.
We use $E=50000$ and $P=2500$ for this problem which implies that we use $Q=20$ gradient descent steps to solve the subproblems in penalty method (${\rm P}^\infty$) and $Q^A=Q^B=10$ gradient descent steps to solve the subproblems in augmented Lagrangian method (${\rm AL}^\infty_\infty$). Figure~\ref{fig:gs-error} shows the various errors and run times for this problem. Both methods produce qualitatively similar error curves. The Lagrange multiplier in this case can be thought of as a function of a single variable since it is a member of $W^{1,2}(\partial\Omega;\mathbb R)$ and $\partial\Omega$ is homeomorphic to a closed curve in $\mathbb R^2$. But for a convenient implementation we represent it as a function of two variables (or on $\Omega$) as seen in table~\ref{tab:network}. This does not cause any practical issues since we never encounter the multiplier anywhere expect $\partial \Omega$ during the run time of ${\rm AL}^\infty_\infty$ i.e. we only optimize the multiplier on $\partial\Omega$. According to table~\ref{tab:network} the size of the multiplier network $b$ is significantly smaller than the size of the solution network $a$ in this case. These choices reflect the fact that inherently the solution and the multiplier are functions of two and one variables respectively. Even though on a machine, \eqref{eq:seq-uncon-al-net} is more expensive to solve than \eqref{eq:seq-uncon-net}, since \eqref{eq:xi-update} is much cheaper to solve than \eqref{eq:seq-uncon-net}, the augmented Lagrangian is ultimately much faster than its counterpart in this case as is seen in figure~\ref{fig:gs-error}. 
\subsection{Beltrami fields}
\begin{figure}[!ht]
    \centering
\includegraphics[scale=0.32]{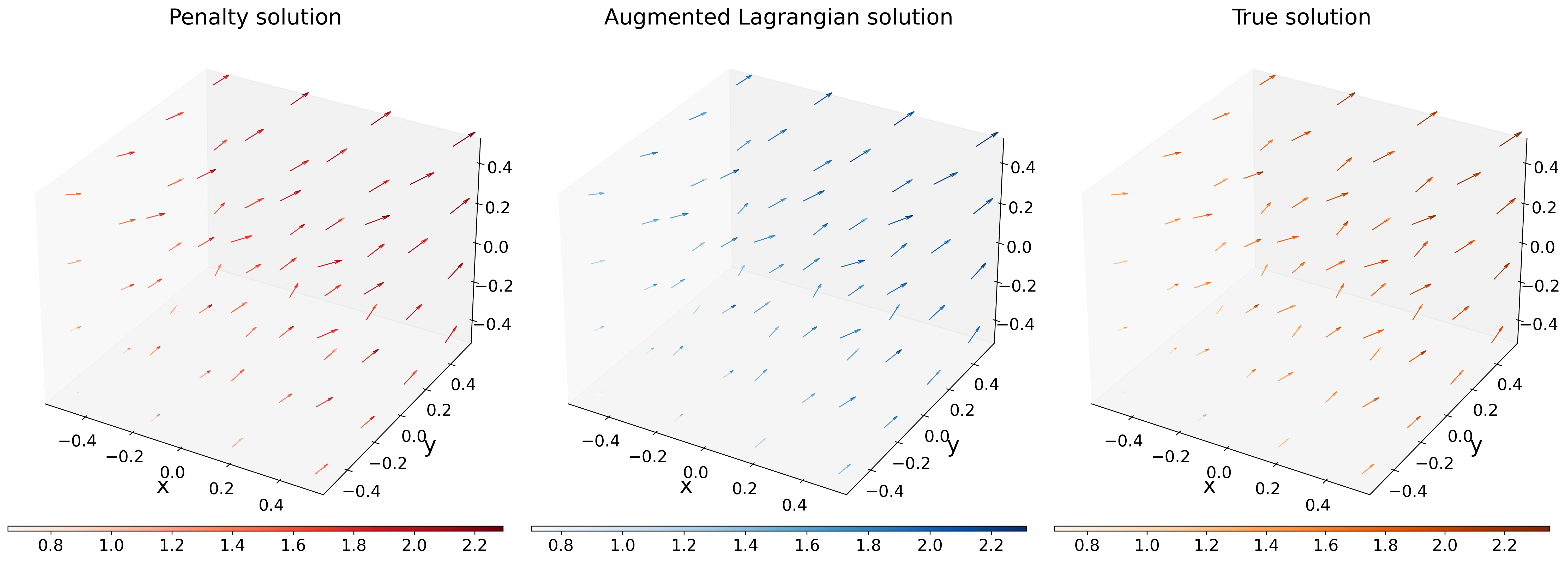}
    \caption{Solutions to the Beltrami field problem. The arrows-lengths have been normalized for visual clarity. The colorbars depict the magnitude of the vectors.}
    \label{fig:bel}
\end{figure}

\begin{figure}[!ht]
    \centering
\includegraphics[scale=0.32]{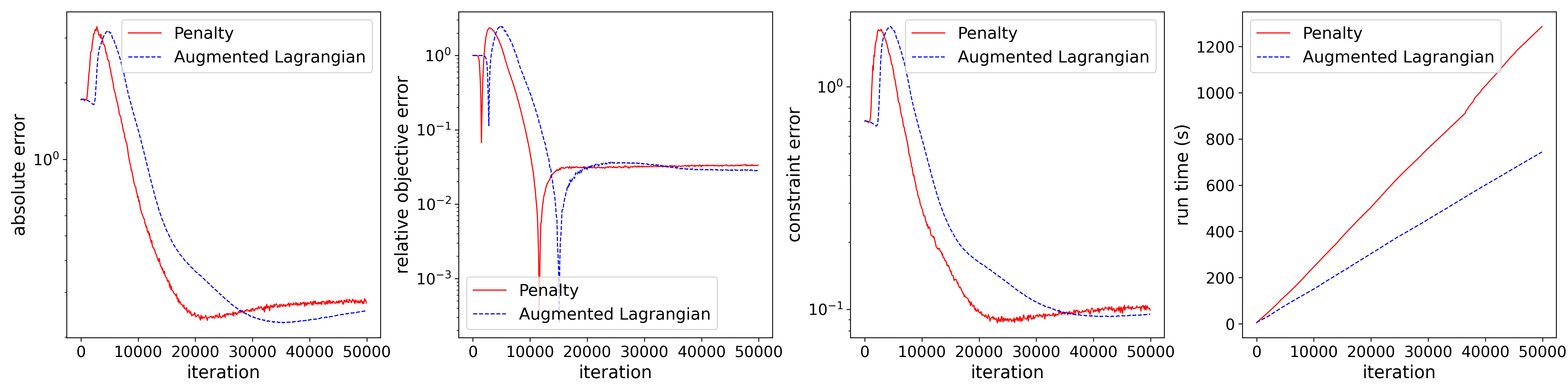}
    \caption{Errors and run times for the Beltrami field problem as functions of gradient descent steps. The errors have been plotted in a semilog fashion. All quantities have been plotted every 100 steps.}
    \label{fig:bel-error}
\end{figure}

Figure~\ref{fig:bel} shows the approximate and true solutions for the Beltrami field problem.
 We use $E=50000$ and $P=2500$ for this problem which implies that we use $Q=20$ gradient descent steps to solve the subproblems in penalty method (${\rm P}^\infty$) and $Q^A=Q^B=10$ gradient descent steps to solve the subproblems in augmented Lagrangian method (${\rm AL}^\infty_\infty$). Figure~\ref{fig:bel-error} shows the various errors and run times for this problem. In order to enforce Gauss's law we represent the solution magnetic field $\hat u$ as the curl of a vector potential and we represent this vector potential as a neural network $H^{\mathcal A}_{\eta}$, 
 \begin{align}
     \hat u = \nabla\times H^{\mathcal A}_{\eta}
 \end{align}
 This clearly allows many such vector potentials to generate solutions to our problem but we do not concern ourselves with gauge-fixing since we are only interested in the magnetic field $\nabla\times H^{\mathcal A}_{\eta}$ rather than the potential $H^{\mathcal A}_{\eta}$ itself. According to table~\ref{tab:network} the vector potential network is much larger than the Lagrange multiplier network in this case. These choices reflect the fact that the solution is a function defined on a volume while the Lagrange multiplier is a function defined on a surface. Just like the previous problem we implement the multiplier as a function of 3 variables while optimizing its values only the boundary $\partial\Omega$ without causing any practical issues. The functional $f$ in this problem requires integration on a 3D volume. To do so we use a Monte-Carlo estimate with only $10^3$ points to maintain a relatively low computational budget. We do achieve qualitatively decent approximations as seen in figure~\ref{fig:bel} but our computational parsimony along with the use of the curl operator leading to more floating point operations and errors result in higher absolute and constraint errors compared to the other problems discussed here. The relative computational ease of solving \eqref{eq:xi-update} when compared to \eqref{eq:seq-uncon-net} results in a faster performance for the augmented Lagrangian algorithm.   

%% file: var-tex/conclusion.tex
In this work we present some practical implementations of popular constrained optimization algorithms in infinite dimensional Hilbert spaces. Both penalty and augmented Lagrangian methods produce decent, comparable solutions for our toy problems in terms of various metrics. Dimension of the Hilbert space $W$ is an important factor when it comes to the difference in the run times of penalty and augmented Lagrangian algorithms. When $W$ is infinite dimensional we might be able to achieve considerably lower run times for the augmented Lagrangian algorithm compared to the penalty method since updating the multiplier is generally less expensive than solving the subproblem in the penalty method. Some constraints like Gauss's law can be implemented reasonably well through architecture. Different update rules for the Lagrange multiplier lead to different variants of algorithm~\ref{algo:dl-al-infinite} exploring which is a worthwhile topic for future research. The geometry of the objective function in the subproblems and the distribution of their optima also deserve further exploration.     

%% file: var-tex/acknowledge.tex
This work was supported by the Department of Atomic Energy, Government of India, under project no. RTI4001 and the Infosys-TIFR Leading Edge travel grant (2023) and the International Research Connections Fund (2022). The author would like to thank Zhisong Qu, Matthew Hole and Robert Dewar for useful discussions on these topics.

%% file: var-tex/appendix.tex
\subsection{Architecture} \label{ssec-arch}
We use two different architectures here which we refer to as LSTM (long short-term memory) and FF (feed-forward). LSTM type networks have been employed to solve partial differential equations \cite{sirignano2018dgm}, \cite{mandal2023learning} and are useful for avoiding vanishing gradients in deep networks \cite{sherstinsky2020fundamentals}, \cite{vennerod2021long}. We use the same LSTM architecture that appears in \cite{mandal2023learning}. For the sake of completeness we describe this architecture in detail below. Here $\mathbf 0_k$ implies a zero vector of dimension $k$, $\odot$ implies the Hadamard product, $d_I, d_O$ denote the input and output dimensions and $\theta$ represents the ordered set of training parameters. The architecture has two numeric hyperparameters $m, L$ describing the size of individual layers and the number of LSTM blocks respectively. Activation $\mathtt A$ is the elementwise $\tanh$ function and finally, $n^{{\rm LSTM}(m, L, d_I, d_O)}_\theta$ represents the complete network.

\begin{align}
    i\in\{1,2,&\cdots, L\}\\
    \mathtt c_0(\mathbf x)\stackrel{\rm def}{=}&\,\mathbf 0_m\label{eq:layer-c0}\\
    \mathtt h_0(\mathbf x)\stackrel{\rm def}{=}&\,\mathbf 0_{d_I}\label{eq:layer-h0}\\
    \mathtt f_i(\mathbf x) \stackrel{\rm def}{=}& \mathtt A(\mathtt W_f^{(i)}\mathbf x + \mathtt U_f^{(i)}\mathtt h_{i-1}(\mathbf x) + \mathtt b_f^{(i)})\label{eq:layer-f}\\
    \mathtt g_i(\mathbf x) \stackrel{\rm def}{=}& \mathtt A(\mathtt W_g^{(i)}\mathbf x + \mathtt U_g^{(i)}\mathtt h_{i-1}(\mathbf x) + \mathtt b_g^{(i)})\label{eq:layer-g}\\
    \mathtt r_i(\mathbf x) \stackrel{\rm def}{=}& \mathtt A(\mathtt W_r^{(i)}\mathbf x + \mathtt U_r^{(i)}\mathtt h_{i-1}(\mathbf x) + \mathtt b_r^{(i)})\label{eq:layer-r}\\
    \mathtt s_i(\mathbf x) \stackrel{\rm def}{=}& \mathtt A(\mathtt W_s^{(i)}\mathbf x + \mathtt U_s^{(i)}\mathtt h_{i-1}(\mathbf x) + \mathtt b_s^{(i)})\label{eq:layer-s}\\
    \mathtt c_i(\mathbf x) \stackrel{\rm def}{=}&  \mathtt f_i(\mathbf x)\odot \mathtt c_{i-1}(\mathbf x) + \mathtt g_i(\mathbf x)\odot s_i(\mathbf x)\label{eq:layer-c}\\
    \mathtt h_i(\mathbf x) \stackrel{\rm def}{=}& \mathtt r_i(\mathbf x)\odot \mathtt A(\mathtt c_i(\mathbf x))\label{eq:layer-h}\\
    \mathtt d_L(\mathbf x)\stackrel{\rm def}{=}&\mathtt W^\top\mathbf x + \mathtt b\label{eq:layer-final}\\
    n^{{\rm LSTM}(m, L, d_I, d_O)}_\theta \stackrel{\rm def}{=}& \mathtt d_L\circ \mathtt h_L\label{eq:def-LSTM} 
\end{align}
Here $\{\mathtt f_i, \mathtt g_i, \mathtt r_i, \mathtt s_i, \mathtt c_i, \mathtt h_i: i=1,\cdots,L\}\cup\{\mathtt d_L\}$ are the hidden layers and
\begin{align}
    \theta=\{\mathtt W_f^{(i)}, \mathtt U_f^{(i)}, \mathtt b_f^{(i)}, \mathtt W_g^{(i)}, \mathtt U_g^{(i)}, \mathtt b_g^{(i)}, \mathtt W_r^{(i)}, \mathtt U_r^{(i)}, \mathtt b_r^{(i)}, \mathtt W_s^{(i)}, \mathtt U_s^{(i)}, \mathtt b_s^{(i)}:
i=1,\cdots,L\}\cup\{\mathtt W, \mathtt b\}\label{eq:theta-composition}
\end{align}is the set of the trainable parameters. The dimensions of these parameters are given below,
\begin{align}
   \mathtt W_f^{(i)}, 
   \mathtt W_g^{(i)},  \mathtt W_r^{(i)},  \mathtt W_s^{(i)} \in&\;
   \mathbb R^{m\times d_I}\\
   \mathtt U_f^{(i)},
\mathtt U_g^{(i)},
\mathtt U_r^{(i)},
\mathtt U_s^{(i)}\in&
   \begin{cases}\mathbb R^{m\times d_I},\quad\text{ if }i=1 \\
   \mathbb R^{m\times m},\quad\text{otherwise}
   \end{cases}\\
   \mathtt b_f^{(i)},\mathtt b_g^{(i)},\mathtt b_r^{(i)},\mathtt b_s^{(i)}\in&\;\mathbb R^m\\
   \mathtt W\in\mathbb R^{d_O\times m}, \mathtt b\in&\;\mathbb R^{d_O}
\end{align}
which implies the size of the network or cardinality of $\theta$ is $4m[d_I(L+1)+m(L-1)+1]+d_O(m+1)$.

We use $n^{{\rm FF}(m, L, d_I, d_O)}_\phi$ to represent a simple feed-forward network without any skip connections with $\phi$ being the set of trainable parameters. In this case the hyperparameters $m, L$ denote the size of an individual layer and the number of hidden layers respectively. 

\begin{align}
    i\in\{1,&\cdots, L-1\}\\
    \mathtt f_{0}(\mathbf x) \stackrel{\rm def}{=}& \mathtt A(\mathtt W_f^{(0)}\mathbf x + \mathtt b_f^{(0)})\\
    \mathtt f_{i}(\mathbf x) \stackrel{\rm def}{=}& \mathtt A(\mathtt W_f^{(i)}\mathtt f_{i-1}(\mathbf x) + \mathtt b_f^{(i)})\\
    \mathtt d_L(\mathbf x)\stackrel{\rm def}{=}&\mathtt W^\top\mathbf x + \mathtt b\\
    n^{{\rm FF}(m, L, d_I, d_O)}_\phi \stackrel{\rm def}{=}& \mathtt d_L\circ \mathtt f_L\label{eq:def-FF} 
\end{align}
So $\phi$ is given by,
\begin{align}
    \phi=\{\mathtt W^{i}_f, \mathtt b^{i}_f:i=0,1,\cdots L-1\}\cup\{\mathtt W, \mathtt b\}
\end{align}
and the dimensions of these trainable parameters are given below,
\begin{align}
   \mathtt W_f^{(i)} \in&\;\begin{cases}
   &\mathbb R^{m\times d_I},\quad\text{ if }i=0\\
   &\mathbb R^{m\times m},\quad\text{otherwise }\end{cases}\\
   \mathtt b_f^{(i)}\in&\;\mathbb R^m\\
   \mathtt W\in\mathbb R^{d_O\times m}, \mathtt b\in&\;\mathbb R^{d_O}
\end{align}
The size of $n^{{\rm FF}(m, L, d_I, d_O)}_\phi$ is therefore $m[d_I + (L+1) + m(L-1)]+(d_O-1)(m+1)$. 
We list the network architectures and sizes used in our experiments in table~\ref{tab:network}.
\begin{table}[!ht]
\begin{center}
\begin{tabular}{ |c|c|c|c|c|c| }   \hline
    Problem & Algorithm  &$\mathcal A$ &  $a$&  $\mathcal B$&  $b$ \\ 
    \hline
    Minimal surface&${\rm P}^\infty$ &  FF(50, 3, 2, 1) & $5300$  &-  &  -\\
    \hline
    Minimal surface &${\rm AL}^\infty_\infty$ &  FF(50, 3, 2, 1) & $5300$  &FF(50, 3, 1, 1)  &  $5250$\\
    \hline
    Geodesic&${\rm P}^\infty$ &  LSTM(50, 3, 1, 1) & $21051$  &- &  - \\
    \hline
    Geodesic &${\rm AL}^\infty_{\rm F}$ &  LSTM(50, 3, 1, 1) & $21051$  & scalar  &  1\\
    \hline
   Grad-Shafranov &${\rm P}^\infty$ &  LSTM(50, 3, 2, 1) & $21851$  &-&  -\\
    \hline
    Grad-Shafranov & ${\rm AL}^\infty_\infty$ &  LSTM(50, 3, 2, 1) & $21851$  & FF(50, 3, 2, 1) &  $5300$\\
    \hline
    Beltrami field& ${\rm P}^\infty$ &  LSTM(50, 3, 3, 3) & $22753$  &- &  -\\
    \hline
    Beltrami field &${\rm AL}^\infty_\infty$ &  LSTM(50, 3, 3, 3) & $22753$  &FF(50, 3, 3, 3)  &  $5452$\\
    \hline
\end{tabular}
\caption{Networks used in various experiments}
\label{tab:network}
\end{center}
\end{table}
\subsection{Penalty factor}\label{ssec-mu} Recall that we use a stopped geometric sequence as our $\mu_k$ \eqref{eq:mu}. We list the hyperparameters that determine $\mu_k$ in table~\ref{tab:mu}.
\begin{table}[!ht]
    \begin{center}
    \begin{tabular}{ |c|c|c|c|c|c| } 
     \hline
     Problem & Algorithm& $\mu_1$ & $\mu_{\max}$ &  $r$\\ 
     \hline
     Minimal surface& ${\rm P}^\infty$ & $100$ & $5000$  &$1.01$  \\
     \hline
     Minimal surface &${\rm AL}^\infty_\infty$ & $100$ & $5000$  &$1.01$  \\
     \hline
     Geodesic&${\rm P}^\infty$ & $100$ & $500$  &$1.01$  \\
     \hline
     Geodesic &${\rm AL}^\infty_{\rm F}$ & $100$ & $500$  &$1.01$  \\
     \hline
    Grad-Shafranov&${\rm P}^\infty$ & $100$ & $1000$  &$1.01$\\
     \hline
     Grad-Shafranov &${\rm AL}^\infty_\infty$ & $100$ & $1000$  &$1.01$ \\
     \hline
     Beltrami field&${\rm P}^\infty$ & $100$ & $5000$  &$1.01$ \\
     \hline
     Beltrami field &${\rm AL}^\infty_\infty$ & $100$ & $5000$  &$1.01$  \\
     \hline
\end{tabular}
\caption{Hyperparameters of the penalty factor for various experiments}
\label{tab:mu}
\end{center}
\end{table}

\subsection{Learning rate}\label{ssec-rate}
The learning rate $\delta$ depends on 7 hyperparameters which are the initial learning rate ($L_0$), initial decay rate ($D_0$), initial decay steps ($S_0$), tipping point ($T$), final learning rate ($L_1$), final decay rate ($D_1$), final decay steps ($S_1$). We define $\delta$ as,
\begin{align}
    \delta(t)=\begin{cases}
         L_0D_0^{\frac{t\,{\rm mod }\,S_0}{S_0}},\;&t<T\\
         L_1D_1^{\frac{t-T}{S_1}},\;&t\ge T\\
    \end{cases}
\end{align}
We list these parameters for our experiments in table~\ref{tab:rate}. We use,
\begin{align}
    &S_0 = \frac{2E}{P}\\
    &T=\left\lfloor\frac{2E(\mu_{\max}-\mu_1)}{Pr}\right\rfloor\\
    &S_1=E-T
\end{align}
For definitions of $E,P$ see section~\ref{sec-appendix}. In case $T>E$, we never reach the tipping point and hence do not list $L_1, D_1$.
\begin{table}[!ht]
\begin{center}
\begin{tabular}{ |c|c|c|c|c|c| } 
    \hline
    Problem & Algorithm& $L_0$ & $D_0$ &  $L_1$ & $D_1$ \\ 
    \hline
    Minimal surface& ${\rm P}^\infty$ & $10^{-4}$ & $10^{-1}$  &- &- \\
    \hline
    Minimal surface &${\rm AL}_\infty^\infty$ & $10^{-4}$ & $2\times10^{-1}$  &- &- \\
    \hline
    Geodesic&${\rm P}^\infty$ & $10^{-3}$ & $10^{-1}$  &$10^{-4}$ &$10^{-2}$ \\
    \hline
    Geodesic &${\rm AL}^\infty_{\rm F}$ & $10^{-3}$ & $10^{-1}$  &$10^{-4}$ &$10^{-2}$ \\
    \hline
   Grad-Shafranov&${\rm P}^\infty$ & $10^{-4}$ & $10^{-1}$  &$10^{-6}$ &$10^{-2}$ \\
    \hline
    Grad-Shafranov &${\rm AL}_\infty^\infty$ & $10^{-4}$ & $10^{-1}$  &$10^{-6}$ &$10^{-2}$ \\
    \hline
    Beltrami field&${\rm P}^\infty$ & $10^{-4}$ & $10^{-1}$  &- &- \\
    \hline
    Beltrami field &${\rm AL}_\infty^\infty$ & $10^{-4}$ & $10^{-1}$  &- &- \\
    \hline
   \end{tabular}
   \caption{ Hyperparameters of the learning rate for various experiments}
   \label{tab:rate}
\end{center}
\end{table}